\documentclass[notitlepage,11pt]{article}
\usepackage{amssymb}
\catcode`\@=11
\@addtoreset{equation}{section}

\catcode`\@=12

\usepackage{latexsym}
\usepackage{amsmath}
\usepackage{amsfonts}
\usepackage{amssymb}
\usepackage{mathrsfs}
\usepackage{comment}
\usepackage{tikz}

\def\R{\mathbb{R}}
\def\S{\mathbb{S}}
\def\f{\varphi}
\def\eps{\varepsilon}

\def\proof{\noindent{\textbf{Proof. }}}
\def\QED{\hfill {$\square$}\goodbreak \medskip}

\newtheorem{Theorem}{Theorem}[section]
\newtheorem{Lemma}[Theorem]{Lemma}
\newtheorem{Proposition}[Theorem]{Proposition}
\newtheorem{Corollary}[Theorem]{Corollary}
\newtheorem{Remark}[Theorem]{Remark}

\begin{document}

\title{Rellich inequalities with weights}

\author{Paolo Caldiroli \and Roberta Musina}

\date{}

\maketitle

\begin{abstract}
Let $\Omega$ be a cone in $\R^{n}$ with $n\ge 2$. For every fixed $\alpha\in\R$ we find the best constant in the Rellich inequality $\int_{\Omega}|x|^{\alpha}|\Delta u|^{2}dx\ge C\int_{\Omega}|x|^{\alpha-4}|u|^{2}dx$ for $u\in C^{2}_{c}(\overline\Omega\setminus\{0\})$. We also estimate the best constant for the same inequality on $C^{2}_{c}(\Omega)$. Moreover we show improved Rellich inequalities with remainder terms involving logarithmic weights on cone-like domains.
\\
\\
\textit{Keywords:} {Rellich inequality, weight, Laplace-Beltrami operator}\\
{\it 2010 Mathematics Subject Classification:} 26D10, 47F05
\end{abstract}

\section*{Introduction}

Let $\alpha$ be a fixed real number and let $\Omega$ be a domain in $\R^{n}$, $n\ge 2$, with $0\in\partial\Omega$. If $\Omega\neq
\R^n\setminus\{0\}$ we assume that 
$\partial\Omega\setminus\{0\}$ is regular enough.
We deal with a class of inequalities of the form
\begin{equation}
\label{eq:intro-inequality}
\int_{\Omega}|x|^{\alpha}|\Delta u|^{2}dx\ge C\int_{\Omega}|x|^{\alpha-4}|u|^{2}dx
\end{equation}
where $u\in C^2(\overline{\Omega})$ runs into suitable classes of functions vanishing on $\partial\Omega$. 
More precisely, we consider the following cases:
\begin{itemize}
\item[(i)] $u$ vanishes on $\partial\Omega$ and in a neighborhood of $0$ and of $\infty$. We will denote this functional space as $C^{2}_{c}(\overline\Omega\setminus\{0\})$ and we will refer to this situation as the \emph{Navier case}. 
\item[(ii)] $u$ has compact support in $\Omega$.  We will denote this functional space as $C^{2}_{c}(\Omega)$ and we will refer to this situation as 
the \emph{Dirichlet case}.
\end{itemize}

Our first goal is to evaluate the best constants 
\begin{equation}
\label{eq:Navier-constant}
\begin{split}
\mu_{N}(\Omega;\alpha):=&\inf_{\scriptstyle u\in C^{2}_{c}(\overline\Omega\setminus\{0\})
\atop\scriptstyle u\ne 0}\frac{\displaystyle\int_{\Omega}|x|^{\alpha}|\Delta u|^{2}dx}{\displaystyle\int_{\Omega}|x|^{\alpha-4}|u|^{2}dx}~,\\
\mu_{D}(\Omega;\alpha):=&\inf_{\scriptstyle u\in C^{2}_{c}(\Omega)
\atop\scriptstyle u\ne 0}\frac{\displaystyle\int_{\Omega}|x|^{\alpha}|\Delta u|^{2}dx}{\displaystyle\int_{\Omega}|x|^{\alpha-4}|u|^{2}dx}~.
\end{split}
\end{equation}
Since the exponents in the power-type weights make the inequality (\ref{eq:intro-inequality}) and the ratios in (\ref{eq:Navier-constant}) 
invariant under dilations, we focus our attention to the case of dilation invariant domains. More precisely, we study (\ref{eq:intro-inequality}) on \emph{cones} $\Omega=\mathcal{C}_{\Sigma}$, where $\Sigma$ is a domain in the unit sphere $\S^{n-1}$ and
\begin{equation}
\label{eq:CSigma}
\mathcal{C}_{\Sigma}:=\left\{~r\sigma~|~r>0~,~\sigma\in\Sigma~\right\}.
\end{equation}
As special cases we include $\Omega=\R^{n}\setminus\{0\}$ (corresponding to $\Sigma=\S^{n-1}$) and $\Omega=\R^{n}_{+}$, a homogeneous half-space (corresponding to $\Sigma=\S^{n-1}_{+}$, a half-sphere). 

We provide an explicit formula for the best constant in the Navier case and estimates from below in the Dirichlet case (see Theorem \ref{T:Rellich-cone} and Corollary \ref{C:Rellich-cone}). 

We also consider {\em cone-like domains}, that are of the following kind:
\begin{itemize}
\item[$i)$]
Bounded domains of the form $\Omega=\mathcal{C}_{\Sigma}\cap\mathbb{B}^{n}$, where $\mathbb{B}^{n}$ denotes the unit ball in $\R^{n}$. This class of domains includes the punctured ball $\mathbb{B}^{n}\setminus\{0\}$, and the half-ball $\mathbb{B}^{n}_{+}=\mathbb{B}^{n}\cap\R^{n}_{+}$.
\item[$ii)$]
Exterior domains of the form $\Omega=\mathcal{C}_{\Sigma}\setminus\overline{\mathbb{B}^{n}}$. 
\end{itemize}

For these domains we will prove inequalities with remainder terms involving logarithmic weights, both in the Navier and in the Dirichlet cases, with optimal constants (see Theorems \ref{T:log} and \ref{T:log-sharp}).

Inequalities of the form (\ref{eq:intro-inequality}) are known in the literature as Rellich-type inequalities. 
Even if they are less studied than the corresponding lower order inequalities, nowadays constitute a fecund field of research. 
The prototype case is the inequality:
\begin{equation}
\label{eq:classical-Rellich}
\int_{\R^{n}}|\Delta u|^2~dx\ge \mu_{n}\int_{\R^{n}}|x|^{-4}|u|^2~dx\quad\textrm{for any $u\in C^{2}_c({\R^{n}}\setminus\{0\})$.}
\end{equation}
It was proved by Rellich in 1953 (see \cite{Rel54} and the posthumous paper \cite{Rel69}) with the optimal constant
\begin{equation}
\label{eq:RN}
\mu_{n}=\mu_{D}(\R^{n}\setminus\{0\};0)=\mu_{N}(\R^{n}\setminus\{0\};0)=\begin{cases}
\left(\frac{n(n-4)}{4}\right)^{2}&\displaystyle\textrm{if $n\ne 2$}\\
0&\textrm{if $n=2$.}
\end{cases}
\end{equation}
Subsequently many authors (see \cite{All}, \cite{Ben}, \cite{DavHin}, \cite{EvaLew}) studied the more general version
\begin{equation}
\label{eq:alpha-Rellich}
\int_{\R^{n}}|x|^{\alpha}|\Delta u|^2~dx\ge \mu_{n,\alpha}\int_{\R^{n}}|x|^{\alpha-4}|u|^2~dx\quad\textrm{for any $u\in C^{2}_c({\R^{n}}\setminus\{0\})$,}
\end{equation}
under some restrictions on $\alpha$ (see also \cite{BerFarFerGaz}, \cite{GazGruMit}, \cite{GazGruSwe}
and \cite{Mit} for related results).

As a corollary to Theorem \ref{T:Rellich-cone}, we provide a complete answer on the validity of (\ref{eq:alpha-Rellich}) for all $\alpha\in\R$ and in any 
dimension $n\ge 2$.
Moreover we compute the sharp value of the best constant
$\mu_{n,\alpha}:=\mu_{D}(\R^{n}\setminus\{0\};\alpha)=\mu_{N}(\R^{n}\setminus\{0\};\alpha)$, that is
\begin{equation}
\label{eq:CNalpha}
\mu_{n,\alpha}
=\min_{k\in\mathbb{N}\cup\{0\}}\left|\gamma_{n,\alpha}+k(n-2+k)\right|^{2}~\!,
\end{equation}
where
\begin{equation}
\label{eq:gammaNalpha}
\gamma_{n,\alpha}=\left(\frac{n-2}{2}\right)^{2}-\left(\frac{\alpha-2}{2}\right)^{2}.
\end{equation}
Thus for every dimension $n\ge 2$ there is an unbounded sequence $(\alpha_{k})_{k\in\mathbb Z}$ of integers such that the Rellich inequality (\ref{eq:alpha-Rellich}) holds with a positive constant $\mu_{n,\alpha}$ if and only if $\alpha\ne\alpha_{k}$ for any $k\in\mathbb Z$. Differently from the best constant in the Hardy inequality for the $L^2$ norm of $|x|^{\alpha/2}|\nabla u|$ (see for instance \cite{Mit}),
the function $\alpha\mapsto\mu_{n,\alpha}$ continues to oscillate as $\alpha\to\pm\infty$.
\bigskip

\centerline{
 \begin{tikzpicture}[xscale=0.7,yscale=0.7]
    \draw[->] (-7.2,0) -- (7.2,0) node[right] {\footnotesize$\alpha$};
    \draw[->] (-2,-0.2) -- (-2,7.6); 
\draw[thick] (0,0) parabola (1.4,0.55) parabola bend (2,0) (3.1,2.2) parabola bend (4,0) (5.1,6.1) parabola bend (6,0) (6.8,7.4);
\draw[thick] (0,0) parabola (-1.4,0.55) parabola bend (-2,0) (-3.1,2.2) parabola bend (-4,0) (-5.1,6.1) parabola bend (-6,0) (-6.8,7.4);
\filldraw(-2,0) circle (1pt) node [below right] {\footnotesize{$0$}};
\filldraw(-4,0) circle (1pt) node [below] {\footnotesize{$-2$}};
\filldraw(-6,0) circle (1pt) node [below] {\footnotesize{$-4$}};
\filldraw(0,0) circle (1pt) node [below] {\footnotesize{$2$}};
\filldraw(2,0) circle (1pt) node [below] {\footnotesize{$4$}};
\filldraw(4,0) circle (1pt) node [below] {\footnotesize{$6$}};
\filldraw(6,0) circle (1pt) node [below] {\footnotesize{$8$}};
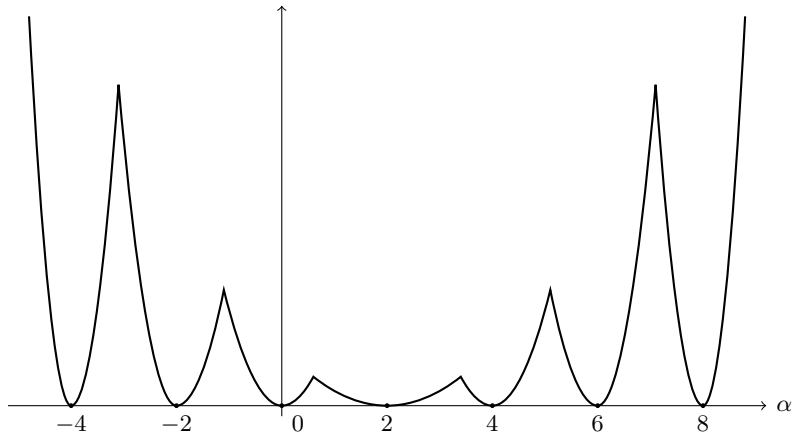
\draw(0,-1) node [below] {\footnotesize{\emph{Graph of $\alpha\mapsto \mu_{2,\alpha}$
in dimension $n=2$}}};
\end{tikzpicture}}
\medskip

Notice that the set $\{k(n-2+k)~|~k=0,1,2,...\}$ is the spectrum of the Laplace-Beltrami operator on the sphere. Thus (\ref{eq:CNalpha}) shows that a resonance phenomenon occurs, that is, the Rellich inequality fails on $\R^n$ if and only if $-\gamma_{n,\alpha}$ is an eigenvalue of $-\Delta_{\sigma}$ on $H^1(\S^{n-1})$.

We point out that in Rellich-type inequalities on the whole space,
the radial and the angular part of $\Delta u$ have independent roles.
Actually no symmetrization or rearrangement argument can
be used to study the minimizations problems in (\ref{eq:Navier-constant}).
In fact, we can show that considering just \emph{radial} functions we have
$$
\inf_{\scriptstyle u\in C^{2}_c({\R^{n}}\setminus\{0\})\atop\scriptstyle u\ne 0,~u(x)=u(|x|)}
\frac{\displaystyle\int_{\R^{n}}|x|^{\alpha}|\Delta u|^{2}dx}{\displaystyle\int_{\R^{n}}|x|^{\alpha-4}|u|^{2}dx}=\gamma_{n,\alpha}^{2}
$$
whereas, taking \emph{nonradial} functions, we have
$$
\inf_{\scriptstyle u\in C^{2}_c({\R^{n}}\setminus\{0\})\atop\scriptstyle u\ne 0,~\int_{\partial B_{r}}u=0~\forall r>0}
\frac{\displaystyle\int_{\R^{n}}|x|^{\alpha}|\Delta u|^2~dx}{\displaystyle\int_{\R^{n}}|x|^{\alpha-4}|u|^2~dx}=\min_{k\in\mathbb{N}}\left|\gamma_{n,\alpha}+k(n-2+k)\right|^{2}.
$$
For instance, in case $\alpha=0$ the Rellich inequality holds for all mappings in $C^{2}_{c}(\R^{n}\setminus\{0\})$ if and only if $n\ne 2,4$ and for all radial mappings 
when $n=2$. In addition, it holds on functions that are orthogonal to radially symmetric maps in any dimension $n\ne 2$.

In case of a cone $\mathcal{C}_{\Sigma}$ as in (\ref{eq:CSigma}), with $\Sigma$ strictly contained in $\S^{n-1}$, 
we have to distinguish between the Navier case and the Dirichlet one. The resonance phenomenon pointed out when
$\Omega=\R^n\setminus\{0\}$
can be observed on any dilation-invariant domain. 
Actually we can prove that 
\begin{equation}
\label{eq:result}
\mu_{N}(\mathcal{C}_{\Sigma};\alpha)=\textrm{dist}\left(-\gamma_{n,\alpha},\Lambda(\Sigma)\right)^{2}~\!,
\end{equation}
where $\Lambda(\Sigma)$ is the Dirichlet spectrum of the Laplace-Beltrami operator on $\Sigma$. In particular
$\mu_{N}(\mathcal{C}_{\Sigma};\alpha)>0$
if and only if $-\gamma_{n,\alpha}$ is not an eigenvalue of $-\Delta_{\sigma}$ on $H^1_0(\Sigma)$.
\bigskip

\centerline{
\begin{tikzpicture}[xscale=1.2,yscale=1.2]
\draw[->] (-0.2,0) -- (4.2,0); 
\draw[->] (0,-0.2) -- (0,4.2); 
\draw[color=black,thick] (1.5708,0) parabola (2.4836,0.72) parabola[bend at end] (3.1416,0);
\draw[color=black,thick,xscale=2,xshift=-0.7854cm]
(0.90000, 4.18689)--(0.92795, 3.47990)--(0.95590, 2.89112)--(0.98385, 2.39966)--(1.01180, 1.98865)--(1.03975, 1.64444)--(1.06770, 1.35589)--(1.09565, 1.11389)--(1.12360, 0.91092)--(1.15155, 0.74080)--(1.17950, 0.59840)--(1.20745, 0.47942)--(1.23540, 0.38031)--(1.26335, 0.29806)--(1.29130, 0.23016)--(1.31925, 0.17448)--(1.34720, 0.12924)--(1.37515, 0.09290)--(1.40310, 0.06417)--(1.43105, 0.04196)--(1.45900, 0.02532)--(1.48695, 0.01345)--(1.51490, 0.00565)--(1.54285, 0.00134)--(1.57080,0.00000);
    \filldraw(1.5708,0) circle(0.8pt) node[below] {\footnotesize{$\frac{\pi}{2}$}};
    \filldraw(3.1416,0) circle(0.8pt) node[below] {\small{$\pi$}};
    \filldraw(2.4836,0) circle(0.8pt) node[below] {\footnotesize{$\sqrt{\frac{5}{8}}\pi$}};
    \filldraw(0,0.72) circle(0.8pt) node[left] {\footnotesize{$\frac{9}{25}$}};
    \draw[thin,dashed](0,0.72)--(2.4836,0.72)--(2.4836,0);
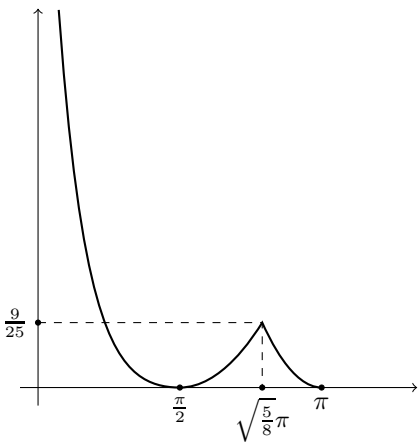
\draw(2,-0.8) node[below]{\footnotesize{Graph of $\theta\mapsto \mu_{N}(\mathcal{C}_{\Sigma_{\theta}};0)$ where $\Sigma_{\theta}$ is an arc in $\S^{1}$ of length $2\theta$.}};
\end{tikzpicture}}

\bigskip

A relevant case is when the cone $\mathcal{C}_{\Sigma}$ is an
homogeneous half-space 
$\R^{n}_{+}$. Computing the Dirichlet spectrum on a half-sphere
(see Proposition \ref{P:spectrum}), from (\ref{eq:result}) we infer
$$
\mu_{N}(\R^{n}_{+};\alpha)=\min_{k\in\mathbb{N}}\left|\gamma_{n,\alpha}+k(n-2+k)\right|^{2}.
$$
In the {\em Dirichlet case} we can prove that
$$
\mu_{D}(\mathcal{C}_{\Sigma};\alpha)>\mu_{N}(\mathcal{C}_{\Sigma};\alpha)
$$
as soon as $\Sigma\ne\S^{n-1}$. In particular,  the Rellich constant on proper cones with Dirichlet
boundary conditions is always positive. For
instance, we can show that
$$
\mu_{D}(\R^{2}_{+};0)>\frac{9}{25}
$$
whereas $\mu_{N}(\R^{2}_{+};0)=0$.

We emphasize the fact that the dependence of the constant $\mu_{N}(\Omega;\alpha)$ with respect to the domain $\Omega$ in general exhibits no monotonicity property. Instead, in the Dirichlet case, one easily sees that $\mu_{D}(\Omega;\alpha)\ge\mu_{D}(\Omega';\alpha)$ if $\Omega\subset\Omega'$.

A more detailed analysis about Rellich inequalities on cones is developed in Section~\ref{S:cones}.
The special cases previously discussed for the whole space or for a half-space are displayed in Section~\ref{S:special}.

Next we deal with cone-like domains.
Here, for the sake of simplicity, we limit ourselves to state some of our results when the domain is either the punctured ball or the complement of the ball. Let $\mu_{n,\alpha}$ be given by (\ref{eq:CNalpha}) and let
\begin{equation}
\label{eq:gammabar}
\overline\gamma_{n,\alpha}:=\left(\frac{n-2}{2}\right)^{2}+\left(\frac{\alpha-2}{2}\right)^{2}.
\end{equation}
If $\Omega=\mathbb{B}^{n}\setminus\{0\}$ or 
$\Omega=\R^{n}\setminus\overline{\mathbb{B}^{n}}$ then the following inequalities hold with sharp constants:
\begin{itemize}
\item[(i)] [Navier case]
for every $u\in C^{2}_{c}(\overline\Omega\setminus\{0\})$ one has
$$
\int_{\Omega}|x|^{\alpha}|\Delta u|^{2}~\!dx-\mu_{n,\alpha}\int_{\Omega}|x|^{\alpha-4}|u|^2~\!dx
\ge\frac{\overline\gamma_{n,\alpha}}{2}\int_{\Omega}|x|^{\alpha-4}|\log|x||^{-2}|u|^2~\!dx~\!;
$$
\item[(ii)] [Dirichlet case]
for every $u\in C^{2}_{c}(\Omega)$ one has
\begin{eqnarray*}
\int_{\Omega}|x|^{\alpha}|\Delta u|^{2}~\!dx-
\mu_{n,\alpha}\int_{\Omega}|x|^{\alpha-4}|u|^2~\!dx
&\ge&
\frac{\overline\gamma_{n,\alpha}}{2}\int_{\Omega}|x|^{\alpha-4}|\log|x||^{-2}|u|^2~\!dx\\
&&+\frac{9}{16}\int_{\Omega}|x|^{\alpha-4}|\log|x||^{-4}|u|^2~\!dx.
\end{eqnarray*}
\end{itemize}
For more general cone-like domains we also consider mixed boundary conditions, precisely of Dirichlet type on the ``radial'' boundary, and of Navier type on the ``angular'' boundary. A wider discussion is contained in Section~\ref{S:log}.

In a forthcoming paper \cite{CM2} we will study second order interpolation inequalities of Caffarelli-Kohn-Nirenberg type \cite{CKN} and reletad noncompact semilinear problems. 


\section{Auxiliary problems on spherical domains}
\label{S:Auxiliary}

For any domain $\Sigma\subseteq\S^{n-1}$ we let $H^1_0(\Sigma)$ to be the Sobolev space of measurable functions on 
$\Sigma$ obtained by completing $C^{2}_{c}(\Sigma)$ with respect to the scalar product
$$
\langle \f, \psi \rangle :=\int_\Sigma \nabla_{\!\sigma} \f\cdot\nabla_{\!\sigma}\psi~d\sigma
+\int_\Sigma \f\psi~d\sigma~\!.
$$
Notice that $H^1_0(\S^{n-1})=H^1(\S^{n-1})$. We denote by $\Lambda(\Sigma)$ the spectrum of the Laplace-Beltrami operator $-\Delta_{\sigma}$ on $H^1_0(\Sigma)$ and by $\lambda_{\Sigma}$ the smallest eigenvalue in $\Lambda(\Sigma)$. Recall that $\lambda_{\Sigma}$ is positive if and only if $\S^{n-1}\setminus\Sigma$ has positive capacity. In this case $\lambda_{\Sigma}$ equals the Poincar\'e constant on $H^1_0(\Sigma)$.

We fix a constant $\gamma\in\R$, and we define the differential operator
$$
L\f=-\Delta_{\sigma} \f+\gamma \f~\!.
$$
In the present section we study the following minimization problems:
$$
m_{N}(\Sigma;\gamma):=\inf_{\genfrac{}{}{0pt}{}{\scriptstyle{\f\in H^2\cap H^1_0(\Sigma)}}{\scriptstyle{\f\neq0}}}~
\frac{\displaystyle{\int_\Sigma|L\f|^2d\sigma}}{\displaystyle\int_\Sigma |\f|^2d\sigma}~,\quad
m_{D}(\Sigma;\gamma):=\inf_{\genfrac{}{}{0pt}{}{\scriptstyle{\f\in H^2_0(\Sigma)}}{\scriptstyle{\f\neq 0}}}~\frac{\displaystyle{\int_\Sigma|L\f|^2d\sigma}}{\displaystyle\int_\Sigma |\f|^2d\sigma}~\!.
$$
Notice that $m_{D}(\Sigma;\gamma)\ge m_{N}(\Sigma;\gamma)$. 
Our first result concerns the lowest infimum $m_{N}(\Sigma;\gamma)$.

\begin{Proposition}
\label{P:auxiliary0}
One has that 
\begin{equation}
\label{eq:mequality}
m_{N}(\Sigma;\gamma)= \mathrm{dist}(-\gamma,\Lambda(\Sigma))^2~\!,
\end{equation}
and $m_{N}(\Sigma;\gamma)$ is always achieved. More precisely,
$\f\in H^{2}\cap H^{1}_{0}(\Sigma)$ attains $m_{N}(\Sigma;\gamma)$ if and only if $\f$ is an eigenfunction relative to the eigenvalue that achieves the minimal distance of $-\gamma$ from $\Lambda(\Sigma)$.
\end{Proposition}

\proof
Let $\lambda\in\Lambda(\Sigma)$ and let $\f$ be an eigenfunction relative to the eigenvalue $\lambda$. 
Since $L\f=(\lambda+\gamma)\f$, then $m_{N}(\Sigma;\gamma)\le(\lambda+\gamma)^2$, and thus $m_{N}(\Sigma;\gamma)\le \textrm{dist}(-\gamma,\Lambda(\Sigma))^2$. 
Therefore it suffices to show that 
\begin{equation}
\label{eq:m2maggiore}
m_{N}(\Sigma;\gamma)\ge\textrm{dist}(-\gamma,\Lambda(\Sigma))^2.
\end{equation} 
If 
$-\gamma$ is an eigenvalue then clearly
$0=m_{N}(\Sigma;\gamma)= \textrm{dist}(-\gamma,\Lambda(\Sigma))^2$
and in addition $m_{N}(\Sigma;\gamma)$ is
achieved by any corresponding eigenfunction. Thus we can assume that
$$
\lambda_{k-1}<-\gamma<\lambda_{k}~\!,
$$
where $\lambda_{k-1}$ and $\lambda_{k}$ are two consecutive eigenvalues
if $-\gamma>\lambda_{\Sigma}$, while $\lambda_{k-1}=-\infty$ if $-\gamma$ is below the spectrum
$\Lambda(\Sigma)$. If $\lambda_{k-1}$ is finite we 
split $H^2\cap H^1_0(\Sigma)$ into the direct sum
$$
H^2\cap H^1_0(\Sigma)=V\oplus V^{\perp}  
$$
where $V$ is the finite-dimensional space spanned by the eigenfunctions relative to the eigenvalues $\lambda<\lambda_k$. Otherwise, we agree that $V=\{0\}$.
Since
$$
\int_\Sigma \f L\f~d\sigma=\int_\Sigma|\nabla_{\sigma}\f|^2~d\sigma+\gamma\int_\Sigma|\f|^2~d\sigma\ge
(\lambda_k+\gamma)\int_{\Sigma}|\f|^2~d\sigma
$$
for any $\f\in V^\perp$, then from the Cauchy-Schwarz inequality
we readily get that
\begin{equation}
\label{eq:m2Nperp}
m^{V^\perp}:=\inf_
{\genfrac{}{}{0pt}{}{\scriptstyle{\f\in V^\perp}}{\scriptstyle{\f\neq 0}}}
~\frac{\displaystyle{\int_\Sigma|L\f|^2d\sigma}}
{\displaystyle\int_\Sigma |\f|^2d\sigma}\ge(\lambda_k+\gamma)^2~\!.
\end{equation}
If $\lambda_{k-1}=-\infty$ then (\ref{eq:mequality}) is proved. If $\lambda_{k-1}$ is finite, namely $V\ne\{0\}$, we show that
\begin{equation}
\label{eq:m2N}
m^{V}:=\inf_
{\genfrac{}{}{0pt}{}{\scriptstyle{\f\in V}}{\scriptstyle{\f\neq 0}}}
~\frac{\displaystyle{\int_\Sigma|L\f|^2d\sigma}}
{\displaystyle\int_\Sigma |\f|^2d\sigma}\ge(\lambda_{k-1}+\gamma)^2~\!.
\end{equation}
Indeed, fix an $L^{2}$-orthonormal basis $\{\f_1,...,\f_h\}$ of $V$, made by eigenfunctions. Any function $\f\in V$ can be written as
$$
\f=\sum_{j=1}^h a_j\f_j
$$
for some $a_{1},...,a_{h}\in\R$. Let $\lambda_{k_j}$ be the eigenvalue relative to $\f_j$. Since $\lambda_{k_{j}}+\gamma\le\lambda_{k-1}+\gamma<0$, we have that
$$
\int_\Sigma|L\f|^2d\sigma=\sum_{j=1}^{h}a_{j}^{2}(\lambda_{k_j}+\gamma)^2\ge(\lambda_{k-1}+\gamma)^{2}\sum_{j=1}^{h}a_{j}^{2}=(\lambda_{k-1}+\gamma)^{2}\int_\Sigma |\f|^2d\sigma
$$
and then (\ref{eq:m2N}) holds.
\\
In order to obtain (\ref{eq:m2maggiore}) we write any nontrivial $\f\in H^2\cap H^1_0(\Sigma)$ as $\f=\f_V+\f_{V^\perp}$, with
$\f_V\in V$ and $\f_{V^\perp}\in V^\perp$. By orthogonality and by (\ref{eq:m2Nperp})--(\ref{eq:m2N}) we get
\begin{eqnarray*}
\frac{\displaystyle\int_\Sigma|L\f|^2~d\sigma}
{\displaystyle\int_\Sigma |\f|^2~d\sigma}&=&
\frac{\displaystyle\int_\Sigma|L\f_V|^2~d\sigma+\displaystyle\int_\Sigma|L\f_{V^\perp}|^2~d\sigma}
{\displaystyle\int_\Sigma |\f_V|^2~d\sigma+\displaystyle\int_\Sigma |\f_{V^\perp}|^2~d\sigma}\\
&\ge&
\frac{m^{V}\displaystyle\int_\Sigma|\f_V|^2~d\sigma+
m^{V^\perp}\displaystyle\int_\Sigma|\f_{V^\perp}|^2~d\sigma}{
\displaystyle\int_\Sigma |\f_V|^2~d\sigma+\displaystyle\int_\Sigma |\f_{V^\perp}|^2~d\sigma}\\
&\ge& \min\left\{(\lambda_{k-1}+\gamma)^2,(\lambda_k+\gamma)^2\right\} 
=\textrm{dist}(-\gamma,\Lambda(\Sigma))^2,
\end{eqnarray*}
as desired. Hence (\ref{eq:mequality}) is proved. The last claim readily
follows from (\ref{eq:mequality}).
\QED

\begin{Proposition}
\label{P:auxiliary 1}
Assume that $\overline{\Sigma}\ne\S^{n-1}$. Then
$$
m_{N}(\Sigma;\gamma)<m_{D}(\Sigma;\gamma)~\!,
$$
and $m_{D}(\Sigma;\gamma)$ is always achieved in $H^2_0(\Sigma)$.
\end{Proposition}

\proof
Clearly $m_{N}(\Sigma;\gamma)\le m_{D}(\Sigma;\gamma)$.
If $-\gamma\notin\Lambda(\Sigma)$ then $m_{N}(\Sigma;\gamma)>0$ by Proposition
\ref{P:auxiliary0} and hence also $m_{D}(\Sigma;\gamma)$ is positive. Thus it is achieved by some $\f\in H^2_0(\Sigma)$, because of the compact embedding of $H^{2}_{0}(\Sigma)$ into $L^{2}(\Sigma)$.
By contradiction, assume that $m_{N}(\Sigma;\gamma)= m_{D}(\Sigma;\gamma)$. Then $\f$ achieves 
$m_{N}(\Sigma;\gamma)$. Thus, by the last assertion in Proposition \ref{P:auxiliary0}, there exists $\lambda\in\Lambda(\Sigma)$ such that $\f\neq 0$ solves
$$
\begin{cases}
-\Delta_{\sigma}\f=\lambda\f&\textrm{in $\Sigma$}\\
\f\in H^2_0(\Sigma),
\end{cases}
$$
contradicting the local strong maximum principle.
If $-\gamma\in\Lambda(\Sigma)$ we take a proper domain $\Sigma'\subset \S^{n-1}$
containing the closure of $\Sigma$ and such that $-\gamma\notin \Sigma'$.
Since $H^2_0(\Sigma)\subset H^2_0(\Sigma')$ then
$$
m_{D}(\Sigma;\gamma)\ge m_{D}(\Sigma';\gamma)\ge m_{N}(\Sigma';\gamma)>0=m_{N}(\Sigma;\gamma)~\!.
$$
Thus $m_{D}(\Sigma;\gamma)>m_{N}(\Sigma;\gamma)$. Finally, $m_{D}(\Sigma;\gamma)$
is achieved since it is positive, via standard arguments.
\QED

\section{Rellich inequalities on cones}
\label{S:cones}

In this section we investigate Rellich inequalities on cones and we evaluate the best Rellich costant in the Navier case. 

We fix a domain $\Sigma$ in $\S^{n-1}$, with $n\ge 2$, and we let $\mathcal{C}_{\Sigma}$ to be the cone in $\R^{n}$ defined by (\ref{eq:CSigma}). Fixing $\alpha\in\R$ we are interested in the infima $\mu_{N}(\mathcal{C}_{\Sigma};\alpha)$ and $\mu_{D}(\mathcal{C}_{\Sigma};\alpha)$ defined as in (\ref{eq:Navier-constant}). 
We start by noticing that
\begin{equation}
\label{eq:alpha4alpha}
\mu_{N}(\mathcal{C}_{\Sigma};\alpha)=\mu_{N}(\mathcal{C}_{\Sigma};4-\alpha)\quad\textrm{and}\quad\mu_{D}(\mathcal{C}_{\Sigma};\alpha)=\mu_{D}(\mathcal{C}_{\Sigma};4-\alpha).
\end{equation}
Indeed, for every $u\colon\overline{\mathcal{C}_{\Sigma}}\setminus\{0\}\to\R$ we define the function $\hat{u}$ on 
$\overline{\mathcal{C}_{\Sigma}}\setminus\{0\}$ by 
$$
\hat{u}(x)=|x|^{2-n}u(|x|^{-2}x).
$$ 
One has that $u\in C^{2}_{c}(\overline{\mathcal{C}_{\Sigma}}\setminus\{0\})$ if and only if $\hat{u}\in C^{2}_{c}(\overline{\mathcal{C}_{\Sigma}}\setminus\{0\})$ and $u\in C^{2}_{c}(\mathcal{C}_{\Sigma})$ if and only if $\hat{u}\in C^{2}_{c}(\mathcal{C}_{\Sigma})$. Moreover
$$
\int_{\mathcal{C}_{\Sigma}}|x|^{\alpha-4}|\hat{u}|^{2}dx=\int_{\mathcal{C}_{\Sigma}}|x|^{-\alpha}|{u}|^{2}dx\quad\textrm{and}\quad
\int_{\mathcal{C}_{\Sigma}}|x|^{\alpha}|\Delta\hat{u}|^{2}dx=\int_{\mathcal{C}_{\Sigma}}|x|^{4-\alpha}|\Delta{u}|^{2}dx.
$$
Hence (\ref{eq:alpha4alpha}) immediately follows.

Before stating our first result, let us recall that by $\Lambda(\Sigma)$ we denote the spectrum of $-\Delta_{\sigma}$ in $H^{1}_{0}(\Sigma)$. Moreover let $\gamma_{n,\alpha}$ be the number defined in (\ref{eq:gammaNalpha}). 

\begin{Theorem}
\label{T:Rellich-cone}
Let $n\ge 2$ and let $\Sigma$ be a domain in $\S^{n-1}$ of class 
$C^{2}$. Then
\begin{equation}
\label{eq:rellich_cono_alpha}
\mu_{N}(\mathcal{C}_{\Sigma};\alpha)=\mathrm{dist}\left(-\gamma_{n,\alpha},\Lambda(\Sigma)\right)^{2}. 
\end{equation}
Moreover, if $\Sigma\neq \S^{n-1}$ then
$$
\mu_{D}(\mathcal{C}_{\Sigma};\alpha)>\mathrm{dist}\left(-\gamma_{n,\alpha},\Lambda(\Sigma)\right)^{2}. 
$$
\end{Theorem}

Theorem \ref{T:Rellich-cone} will be proved in Subsection \ref{SS:proof}. From the monotonicity property of the mapping $\Sigma\mapsto\mu_{D}(\mathcal{C}_{\Sigma};\alpha)$ we have that
\begin{equation}
\label{eq:monotone}
\mu_{D}(\mathcal{C}_{\Sigma};\alpha)\ge\sup\{\mu_{N}(\mathcal{C}_{\Sigma'};\alpha)~|~\Sigma'\supset\Sigma,~\Sigma'~\textrm{of~class~}C^{2}\}.
\end{equation}
Therefore, from Theorem \ref{T:Rellich-cone}, we infer the next result.

\begin{Corollary}
\label{C:Rellich-cone}
Let $n\ge 2$ and let $\Sigma$ be a domain in $\S^{n-1}$
of class $C^{2}$. If $\Sigma\neq\S^{n-1}$ then
$\mu_{D}(\mathcal{C}_{\Sigma};\alpha)>0$.
\end{Corollary}

In our second main result we show that extremal functions do not exist. This is trivial when $\mu_{N}(\mathcal{C}_{\Sigma};\alpha)$ or $\mu_{D}(\mathcal{C}_{\Sigma};\alpha)$ vanish.  When they are positive we need to introduce suitable Sobolev spaces as follows.
In particular, when $\mu_{N}(\mathcal{C}_{\Sigma};\alpha)>0$, we can define a norm on $C^{2}_{c}(\overline{\mathcal{C}_{\Sigma}}\setminus\{0\})$ by setting
\begin{equation}
\label{eq:norm2}
\|u\|_{2,\alpha} =\left(\int_{\mathcal{C}_{\Sigma}}|x|^{\alpha}|\Delta u|^{2}dx\right)^{1/2}~\!.
\end{equation}
The completion of $C^{2}_{c}(\overline{\mathcal{C}_{\Sigma}}\setminus\{0\})$ with respect to this norm will be denoted $\mathcal{N}^{2}(\mathcal{C}_{\Sigma};\alpha)$.

In the same way, when $\mu_{D}(\mathcal{C}_{\Sigma};\alpha)>0$ we introduce the Sobolev space $\mathcal{D}^{2}(\mathcal{C}_{\Sigma};\alpha)$ as the completion of $C^{2}_{c}(\mathcal{C}_{\Sigma})$ with respect to the norm (\ref{eq:norm2}).
 
\begin{Theorem}
\label{T:Rellich-cone-2}
The infima $\mu_{N}(\mathcal{C}_{\Sigma};\alpha)$ and
$\mu_{D}(\mathcal{C}_{\Sigma};\alpha)$ are never attained.
\end{Theorem}

To prove Theorems \ref{T:Rellich-cone} and \ref{T:Rellich-cone-2} we will use a suitable Emden-Fowler transform, that maps functions defined on $\mathcal{C}_{\Sigma}$ into functions on the cylinder
$$
\mathcal{Z}_{\Sigma}:=\{(s,\sigma)\in\R\times\S^{n-1}~|~s\in\R,~\sigma\in\Sigma\}.
$$
This will be done in the next subsections.

\subsection{The Emden-Fowler transform}
\label{S:EF}

To any $u\in C^{2}_{c}(\overline{\mathcal{C}_{\Sigma}}\setminus\{0\})$ we associate a function $w\colon\overline{\mathcal{Z}_{\Sigma}}\to\R$ via
\begin{equation}
\label{eq:uw}
u(x)=|x|^{\frac{4-n-\alpha}{2}}~\!w\left(-\log|x|,\frac{x}{|x|}\right)
\end{equation}
and we define
$$
Tu:=w.
$$
We denote by $C^{2}_{c}(\overline{\mathcal{Z}_{\Sigma}})$ the space of mappings $w\in C^{2}(\overline{\mathcal{Z}_{\Sigma}})$ such that $w(\cdot,\sigma)=0$ for every $\sigma\in\partial\Sigma$ and $w(s,\cdot)=0$ for $|s|$ large enough. In addition we introduce the differential operator
\begin{equation}
\label{eq:Lw-cone}
Lw=-\Delta_{\sigma}w+\gamma_{n,\alpha}w,
\end{equation}
as in Section \ref{S:Auxiliary}, with $\gamma_{n,\alpha}$ defined in (\ref{eq:gammaNalpha}).

\begin{Lemma}
\label{L:EF1}
If $u\in C^{2}_{c}(\overline{\mathcal{C}_{\Sigma}}\setminus\{0\})$ then $Tu\in C^{2}_{c}(\overline{\mathcal{Z}_{\Sigma}})$. If $u\in C^{2}_{c}(\mathcal{C}_{\Sigma})$ then $Tu\in C^{2}_{c}(\mathcal{Z}_{\Sigma})$. Moreover, setting $w=Tu$, one has
\begin{gather}
\label{eq:uvw}
\int_{\mathcal{C}_{\Sigma}}|x|^{\alpha-4}|u|^2 dx= 
\int_{\mathcal Z_\Sigma}|w|^2 dsd\sigma\\
\label{eq:Delta2}
\int_{\mathcal{C}_{\Sigma}}|x|^{\alpha}|\Delta u|^{2}dx=\int_{\mathcal{Z}_{\Sigma}}|Lw|^{2}dsd\sigma+G(w)
\end{gather}
where
\begin{equation}
\label{eq:Delta2-G}
G(w):=\int_{\mathcal{Z}_{\Sigma}}\left(|w_{ss}|^{2}+2|\nabla_{\sigma}w_{s}|^{2}+2\overline\gamma_{\alpha}|w_{s}|^{2}\right).
\end{equation}
\end{Lemma}

\proof
The first two statements and (\ref{eq:uvw})  are trivial. 
Let 
$v:=|x|^{\frac{n-4+\alpha}{2}}u$ and let $w$ be defined as in (\ref{eq:uw}). We compute
$$
\Delta u=|x|^{\frac{4-n-\alpha}{2}}\left[\Delta v+(4-n-\alpha)|x|^{-1}v_{r}-\gamma_{n,\alpha}|x|^{-2}v\right],
$$
where $v_{r}=|x|^{-1}(x\cdot \nabla v)$ denotes the radial derivative of $v$. Now we go from $v$ to $w$, via the transform
$$
v(x)= w\left(-\log|x|,\frac{x}{|x|}\right).
$$
Denoting $w_s$ and $w_{ss}$ the partial derivatives with respect to the real variable of $w$, since
$$
-|x|v_{r}=w_{s}~,\quad|x|^{2}v_{rr}=w_{ss}+w_{s}~,\quad\nabla_{\!\sigma}v=\nabla_{\!\sigma}w~,\quad\Delta_{\sigma}v=\Delta_{\sigma}w~,
$$
we infer that
$$
\Delta u=|x|^{-\frac{n+\alpha}{2}}\left[-Lw+w_{ss}+(\alpha-2)w_{s}\right]~.
$$
Therefore
$$
\int_{\mathcal{C}_{\Sigma}}|x|^{\alpha}|\Delta u|^{2}dx=\int_{\mathcal{Z}_{\Sigma}}|-Lw+w_{ss}+(\alpha-2)w_{s}|^{2}dsd\sigma.
$$
Since
\begin{gather*}
\int_{\mathcal{Z}_{\Sigma}}ww_{s}dsd\sigma=0~,\quad\int_{\mathcal{Z}_{\Sigma}}w_{ss}w_s~dsd\sigma=0~,\\
\int_{\mathcal{Z}_{\Sigma}}(\Delta_{\sigma}w)w_s~dsd\sigma=-\int_{\mathcal{Z}_{\Sigma}}\nabla_{\!\sigma}w\cdot \nabla_{\!\sigma}w_s~dsd\sigma=-\frac{1}{2}\int_{\mathcal{Z}_{\Sigma}}\partial _s|\nabla_{\!\sigma} w|^2~dsd\sigma = 0~\!,
\end{gather*}
we readily obtain (\ref{eq:Delta2})--(\ref{eq:Delta2-G}). 
\QED

\subsection{Proof of Theorem \ref{T:Rellich-cone}}
\label{SS:proof}

Firstly we prove (\ref{eq:rellich_cono_alpha}). By Proposition \ref{P:auxiliary0}, it suffices to show that 
\begin{equation}
\label{eq:equality-Navier}
\mu_{N}(\mathcal{C}_{\Sigma};\alpha)=m_{N}(\Sigma;\gamma_{n,\alpha}).
\end{equation}
Let $T$ be the Emden-Fowler transform. Fix $u\in C^{2}_{c}(\overline{\mathcal{C}_{\Sigma}}\setminus\{0\})$ and put $w=Tu$. Then use (\ref{eq:Delta2}) and (\ref{eq:Delta2-G}). Since $G(w)\ge 0$, by Proposition \ref{P:auxiliary0} and by (\ref{eq:uvw}) we obtain that
\begin{equation*}
\begin{split}
\int_{\mathcal{C}_{\Sigma}}|x|^{\alpha}|\Delta u|^{2}dx&\ge\int_{\mathcal{Z}_{\Sigma}}|Lw|^{2}dsd\sigma\\
&\ge m_{N}(\Sigma;\gamma_{n,\alpha})\int_{\mathcal{Z}_{\Sigma}}|w|^{2}dsd\sigma =m_{N}(\Sigma;\gamma_{n,\alpha})\int_{\mathcal{C}_{\Sigma}}|x|^{\alpha-4}|u|^{2}dx~\!.
\end{split}
\end{equation*}
Hence $\mu_{N}(\mathcal{C}_{\Sigma};\alpha)\ge m_{N}(\Sigma;\gamma_{n,\alpha})$.
In order to prove the opposite inequality we take a function
$w\in C^{2}_{c}(\overline{\mathcal{Z}_{\Sigma}})$ of the form 
$$
w(s,\sigma)=t^{1/2}\psi(ts)\f_{h}(\sigma)~\!,
$$ 
with $t>0$, $\psi\in C^{2}_c(\R)$ and $(\f_{h})\subset C^{2}_{c}(\overline{\Sigma})$ satisfying
$$
\int_{\Sigma}|L\f_{h}|^{2}d\sigma=m_{N}(\Sigma;\gamma_{n,\alpha})+o(1)~,\quad \int_\Sigma|\f_{h}|^2d\sigma=1~,\quad\int_{-\infty}^\infty |\psi|^2ds=1~\!
$$
with $o(1)\to 0$ as $h\to\infty$.
Then 
$$
\int_{\mathcal{C}_{\Sigma}}|x|^{\alpha-4}|u|^2=1
$$
whereas, by (\ref{eq:Delta2})--(\ref{eq:Delta2-G}),
$$
\int_{\mathcal{C}_{\Sigma}}|x|^{\alpha}|\Delta u|^2=\int_{\Sigma}|L\f_{h}|^{2}+t^{4}\int_{-\infty}^{\infty}|\psi''|^{2}+2t^{2}\int_{-\infty}^{\infty}|\psi'|^{2}\left(\overline\gamma_{n,\alpha}+\int_{\Sigma}|\nabla_{\sigma}\f_{h}|^{2}\right).
$$
Taking $t\to 0$ and $h\to\infty$ we immediately obtain that $\mu_{N}(\mathcal{C}_{\Sigma};\alpha)\le m_{N}(\Sigma;\gamma_{n,\alpha})$. Hence (\ref{eq:equality-Navier}) holds true. In the same way one shows that
$$
\mu_{D}(\mathcal{C}_{\Sigma};\alpha)=m_{D}(\Sigma;\gamma_{n,\alpha}).
$$
Then the conclusion follows from Proposition \ref{P:auxiliary 1}.
\QED

\subsection{The functional spaces $\mathcal{N}^{2}(\mathcal{C}_{\Sigma};|x|^{\alpha}dx)$ and $\mathcal{D}^{2}(\mathcal{C}_{\Sigma};|x|^{\alpha}dx)$}
\label{Ss:spaces}

Assume that $\mu_{N}(\mathcal{C}_{\Sigma};\alpha)>0$. Let $\mathcal{N}^{2}(\mathcal{C}_{\Sigma};|x|^{\alpha}dx)$ be the Hilbert space
endowed with the norm (\ref{eq:norm2}).

Using an interpolation argument we endow the space $H^2\cap H^1_0(\mathcal{Z}_{\Sigma})$ with the equivalent norm
\begin{equation}
\label{eq:norma-w}
\|w\|^{2}=\int_{\mathcal{Z}_{\Sigma}}\left(|\Delta_{\sigma} w|^{2}+|w_{ss}|^{2}+|w|^{2}\right)dsd\sigma~\!.
\end{equation}

For every $u\in C^{2}_{c}(\overline{\mathcal{C}_{\Sigma}}\setminus\{0\})$ let $w=Tu$ be the Emden-Fowler transform of $u$. By Lemma \ref{L:EF1}, we have that $Tu\in H^{2}\cap H^{1}_{0}(\mathcal Z_\Sigma)$. 

\begin{Lemma}
\label{L:EF2}
The operator $T\colon C^{2}_{c}(\overline{\mathcal{C}_{\Sigma}}\setminus\{0\})\to H^{2}\cap H^{1}_{0}(\mathcal Z_\Sigma)$ admits a unique continuous extension on $\mathcal{N}^{2}(\mathcal{C}_{\Sigma};|x|^\alpha dx)$ which is an isomorphism between 
the spaces $\mathcal{N}^{2}(\mathcal{C}_{\Sigma};|x|^\alpha dx)$ and $H^{2}\cap H^{1}_{0}(\mathcal Z_\Sigma)$.
Moreover the equalities (\ref{eq:uvw}) and (\ref{eq:Delta2})--(\ref{eq:Delta2-G}) hold true for every function $u$ in $\mathcal{N}^{2}(\mathcal{C}_{\Sigma};|x|^\alpha dx)$.
\end{Lemma}

\proof
Since $\mu_{N}(\mathcal{C}_{\Sigma};\alpha)>0$ an equivalent norm to $\|\cdot\|_{2,\alpha}$ in $\mathcal{N}^{2}(\mathcal{C}_{\Sigma};|x|^\alpha dx)$ is given by
\begin{equation}
\label{eq:norma2alpha}
\|u\| ^{2}:=\int_{\mathcal{C}_{\Sigma}}|x|^{\alpha}|\Delta u|^{2}dx+\int_{\mathcal{C}_{\Sigma}}|x|^{\alpha-4}|u|^{2}dx.
\end{equation}
By density, equalities (\ref{eq:uvw}) and (\ref{eq:Delta2}) hold true for every $u\in\mathcal{N}^{2}(\mathcal{C}_{\Sigma};|x|^\alpha dx)$. Recalling the definitions of the norms $\|u\| $ and $\|w\|$ given in (\ref{eq:norma2alpha}) and in (\ref{eq:norma-w}), respectively, and using also (\ref{eq:uvw}), we have that
$$
\|u\| ^{2}=\|w\|^{2}+2\gamma_{n,\alpha}\int_{\mathcal Z_\Sigma}|\nabla_{\sigma}w|^{2}+2\int_{\mathcal Z_\Sigma}|\nabla_{\sigma}w_{s}|^{2}+\gamma_{n,\alpha}^{2}\int_{\mathcal Z_\Sigma}|w|^{2}+2\overline\gamma_{n,\alpha}\int_{\mathcal Z_\Sigma}|w_{s}|^{2}.
$$
Hence if $\gamma_{n,\alpha}\ge 0$ then $\|u\| \ge\|Tu\|^{2}$ for all $u\in \mathcal{N}^{2}(\mathcal{C}_{\Sigma};|x|^\alpha dx)$ and the conclusion follows. If $\gamma_{n,\alpha}<0$, using (\ref{eq:Delta2})--(\ref{eq:Delta2-G}), we firstly estimate
$$
\|u\| ^{2}\ge\int_{\mathcal Z_\Sigma}|Lw|^{2}+\int_{\mathcal Z_\Sigma}|w_{ss}|^{2}+\int_{\mathcal Z_\Sigma}|w|^{2}.
$$
Then we use the Young inequality $2ab\le\eps^{-1}a^{2}+\eps b^{2}$ ($\eps>0$, $a,b\in\R$) to estimate
$$
\int_{\mathcal Z_\Sigma}|Lw|^{2}\ge\left(1-\eps^{-1}\right)\int_{\mathcal Z_\Sigma}|\Delta_{\sigma}w|^{2}+(1-\eps)\gamma_{n,\alpha}^{2}\int_{\mathcal Z_\Sigma}|w|^{2}.
$$
Hence we infer that
$$
\|u\| ^{2}\ge\min\left\{1-\eps^{-1},1-\gamma_{n,\alpha}^{2}(\eps-1)\right\}\|Tu\|^{2}
$$
and, fixing $\eps$ with $1<\eps<1+\gamma_{n,\alpha}^{-2}$, the conclusion follows as before.
\QED

In a similar way one has:

\begin{Lemma}
\label{L:EF3}
The Emden-Fowler operator $T$ is an isomorphism between the spaces $\mathcal{D}^{2}(\mathcal{C}_{\Sigma};|x|^\alpha dx)$ and $H^{2}_{0}(\mathcal{Z}_{\Sigma})$.
\end{Lemma}

\begin{Remark}
If $\Sigma\neq\S^{n-1}$, then $H^{2}_{0}(\mathcal{Z}_{\Sigma})$ is
properly contained in $H^{2}\cap H^{1}_{0}(\mathcal{Z}_{\Sigma})$.
Thus $\mathcal{D}^{2}(\mathcal{C}_{\Sigma};|x|^\alpha dx)\subset\mathcal{N}^{2}(\mathcal{C}_{\Sigma};|x|^\alpha dx)$ with strict inclusion, that is, 
$C^{2}_{c}(\mathcal{C}_{\Sigma})$ is not dense in $\mathcal{N}^{2}(\mathcal{C}_{\Sigma};|x|^{\alpha}dx)$. 
\end{Remark}

\subsection{Proof of Theorem \ref{T:Rellich-cone-2}}

Assume that $\mu_{N}(\mathcal{C}_{\Sigma};\alpha)$ is attained by some $u\in\mathcal{N}^{2}(\mathcal{C}_{\Sigma},|x|^{\alpha}dx)$, $u\ne 0$. 
By (\ref{eq:uvw}), (\ref{eq:Delta2}) and by Lemma \ref{L:EF2}, we have that
$$
\mu_{N}(\mathcal{C}_{\Sigma};\alpha)
=\inf_{\genfrac{}{}{0pt}{}{\scriptstyle{w\in H^{2}\cap H^{1}_{0}(\mathcal Z_\Sigma)}}{\scriptstyle{w\neq 0}}}\!\frac{\displaystyle{\int_{\mathcal Z_\Sigma}|Lw|^{2}+G(w)}}{\displaystyle\int_{\mathcal Z_\Sigma}|w|^{2}}.
$$ 
Therefore the infimum at the right hand side is attained by $w=Tu$. Notice that $w_{s}\ne 0$, otherwise $w=0$ and then $u=0$, too. For $t>0$ let
$w^{t}(s,\sigma)=w(ts,\sigma)$.
Then $w^t\ne 0$, $w^{t}\in H^{2}\cap H^{1}_{0}(\mathcal Z_\Sigma)$  for all $t\in(0,1)$  and 
\begin{equation*}
\begin{split}
\mu_{N}(\mathcal{C}_{\Sigma};\alpha)&=
~\frac{\displaystyle{\int_{\mathcal Z_\Sigma}|Lw|^{2}+G(w)}}{\displaystyle\int_{\mathcal Z_\Sigma}|w|^{2}}\le
~\frac{\displaystyle{\int_{\mathcal Z_\Sigma}|Lw^{t}|^{2}+G(w^{t})}}
{\displaystyle\int_{\mathcal Z_\Sigma}|w^{t}|^{2}}\\&=
~\frac{\displaystyle{\int_{\mathcal Z_\Sigma}|Lw|^{2}+
t^{4}\int_{\mathcal Z_\Sigma}|w_{ss}|^{2}+2t^{2}\int_{\mathcal Z_\Sigma}\left(|\nabla_{\sigma}w_{s}|^{2}+
\overline\gamma_{n,\alpha}|w_{s}|^{2}\right)}}{\displaystyle\int_{\mathcal Z_\Sigma}|w|^{2}}\\&<
~\frac{\displaystyle{\int_{\mathcal Z_\Sigma}|Lw|^{2}+\int_{\mathcal Z_\Sigma}|w_{ss}|^{2}+
2\int_{\mathcal Z_\Sigma}\left(|\nabla_{\sigma}w_{s}|^{2}+\overline\gamma_{n,\alpha}|w_{s}|^{2}\right)}}{\displaystyle\int_{\mathcal Z_\Sigma}|w|^{2}}=
~\mu_{N}(\mathcal{C}_{\Sigma};\alpha)~\!,
\end{split}
\end{equation*}
a contradiction. A similar argument holds in order to show that $\mu_{D}(\mathcal{C}_{\Sigma};\alpha)$ is not attained in $\mathcal{D}^{2}(\mathcal{C}_{\Sigma},|x|^{\alpha}dx)$.
\QED

\section{Applications of Theorem \ref{T:Rellich-cone}}
\label{S:special}

\subsection{Rellich inequality on the whole space}

Let $\R^{n}_{0}=\R^{n}\setminus\{0\}$ and set
$$
\mu_{n,\alpha}:=\mu_{N}(\R^{n}_{0};\alpha)=\mu_{D}(\R^{n}_{0};\alpha).
$$
In the next theorem we denote by
$C^{2}_{c,rad}(\R^{n}_{0})$ the set of radial functions in $C^{2}_{c}(\R^{n}_{0})$ and by $C^{2}_{c,rad}(\R^{n}_{0})^{\bot}$ the set of functions $u\in C^{2}_{c}(\R^{n}_{0})$ such that $\int_{\partial B_{r}}u~d\sigma=0$ for every $r>0$.

\begin{Theorem}
Let $n\ge 2$. One has that
\begin{equation}
\label{eq:RRN}
\mu_{n,\alpha}=\min_{k\in\mathbb{N}\cup\{0\}}\left|\gamma_{n,\alpha}+k(n-2+k)\right|^{2}.
\end{equation}
Moreover
\begin{equation}
\label{eq:R-radial}
\inf_{\scriptstyle u\in C^{2}_{c,rad}(\R^{n}_{0})\atop\scriptstyle u\ne 0}\frac{\displaystyle\int_{\R^{n}}|x|^{\alpha}|\Delta u|^{2}dx}{\displaystyle\int_{\R^{n}}|x|^{\alpha-4}|u|^{2}dx}=\gamma_{n,\alpha}^{2}
\end{equation}
whereas
\begin{equation}
\label{eq:R-nonradial}
\inf_{\scriptstyle u\in C^{2}_{c,rad}(\R^{n}_{0})^{\bot}\atop\scriptstyle u\ne 0}
\frac{\displaystyle\int_{\R^{n}}|x|^{\alpha}|\Delta u|^2~dx}{\displaystyle\int_{\R^{n}}|x|^{\alpha-4}|u|^2~dx}=\min_{k\in\mathbb{N}}\left|\gamma_{n,\alpha}+k(n-2+k)\right|^{2}.
\end{equation}
\end{Theorem}

\proof
Equality (\ref{eq:RRN}) is an application of Theorem \ref{T:Rellich-cone} and of the fact that the spectrum of the Laplace-Beltrami operator in $\S^{n-1}$ is given by
$$
\Lambda(\S^{n-1})=\{k(n-2+k)~|~k\in\mathbb{N}\cup\{0\}\}.
$$
Equality (\ref{eq:R-radial}) can be immediately obtained via Emden-Fowler transformation. The same holds for (\ref{eq:R-nonradial}) with the further remark that, arguing as in the proof of Theorem \ref{T:Rellich-cone},
$$
\inf_{\scriptstyle u\in C^{2}_{c,rad}(\R^{n}_{0})^{\bot}\atop\scriptstyle u\ne 0}\frac{\displaystyle\int_{\R^{n}}|x|^{\alpha}|\Delta u|^2~dx}{\displaystyle\int_{\R^{n}}|x|^{\alpha-4}|u|^2~dx}=\inf_{\genfrac{}{}{0pt}{}{\scriptstyle{\f\in C^{2}(\S^{n-1})}}{\scriptstyle{\f\neq 0,~\int_{\S^{n-1}}\f=0}}}~\frac{\displaystyle{\int_{\S^{n-1}}|L\f|^2d\sigma}}{\displaystyle\int_{\S^{n-1}}|\f|^2d\sigma}~\!.
\eqno{\square}
$$
\smallskip

\begin{Remark}
In this remark we take $\alpha=0$. Clearly we recover the classical Rellich inequality (\ref{eq:classical-Rellich}) with the best constant $\mu_{n}$ defined in (\ref{eq:RN}). Moreover we also point out the following inequalities, which hold true in any dimension $n\ge 2$:
\begin{gather*}
\int_{\R^{n}}|\Delta u|^{2}dx\ge\left(\frac{n(n-4)}{4}\right)^{2}\int_{\R^{n}}|x|^{-4}|u|^{2}dx~~\forall u\in C^{2}_{c,rad}(\R^{n}_{0})\\
\int_{\R^{n}}|\Delta u|^{2}dx\ge\left(\frac{n^{2}-4}{4}\right)^{2}\int_{\R^{n}}|x|^{-4}|u|^{2}dx~~\forall u\in C^{2}_{c,rad}(\R^{n}_{0})^{\bot}
\end{gather*}
with sharp constants. 
\end{Remark}

\begin{Remark}
If $\alpha>4-n$ then the weight $|x|^{\alpha-4}$ is locally integrable, and a density argument can be used in order to show that
 $$
 \inf_{\genfrac{}{}{0pt}{}{\scriptstyle{u\in C^{2}_{c}(\R^{n})}}{\scriptstyle{u\neq 0}}}
\frac{\displaystyle\int_{\R^{n}}|x|^\alpha|\Delta u|^{2}dx}{\displaystyle\int_{\R^{n}}|x|^{\alpha-4}|u|^{2}dx}=\mu_{n,\alpha}.
$$
\end{Remark}

We finally observe that
$\mu_{n,\alpha}=0$ if and only if $\alpha\in (n+2\mathbb{N}\cup\{0\})\cup(4-n-2\mathbb{N}\cup\{0\})$.

\subsection{Rellich inequalities on half-spaces}

Denote by $\R^{n}_+$ any homogeneous half-space in $\R^{n}$. From Theorem \ref{T:Rellich-cone} and computing the spectrum of the Laplace-Beltrami operator on the half-sphere $\S^{n-1}_{+}$, the next result follows.

\begin{Theorem}
Let $n\ge 2$. One has that
$$
\mu_{D}(\R^{n}_{+};\alpha)>\mu_{N}(\R^{n}_{+};\alpha)=\min_{k\in\mathbb{N}}\left|\gamma_{n,\alpha}+k(n-2+k)\right|^{2}.
$$
\end{Theorem}

The knowledge of the spectrum of the Laplace-Beltrami operator on the half-sphere $\S^{n-1}_{+}$ could be known but, since we did not find it in the literature, we prove it for the sake of completeness.

\begin{Proposition}
\label{P:spectrum}
One has that $\Lambda(\S^{n-1}_{+})=\{k(k+n-2)~|~k\in\mathbb{N}\}$.
\end{Proposition}

\proof
We have to show that $\Lambda(\S^{n-1}_{+})=\Lambda(\S^{n-1})\setminus\{0\}$. The inclusion 
$$\Lambda(\S^{n-1}_{+})\subset\Lambda(\S^{n-1})\setminus\{0\}$$ easily follows via odd extension
of any Dirichlet eigenfunction on the half-sphere.
Now let $\lambda=k(k+n-2)$ for some $k\in\mathbb{N}$. By well known results (see, e.g., \cite{BerGauMaz}), $\lambda$ is an eigenvalue of $-\Delta_{\sigma}$ in $H^{1}(\S^{n-1})$ whose eigenspace is spanned by $m=m(k)$ functions $H_{k,1},\ldots,H_{k,m}$ with the following properties: the $H_{m,i}$'s are orthonormal and the mapping
$$
x\mapsto |x|^{k}H_{k,i}(x/|x|)\quad\forall x\in\R^{n}
$$
is a polynomial of degree $k$ and it is harmonic on $\R^{n}$. It is also known that
$$
m=\left(\begin{array}{c}n+k-1\\ k\end{array}\right)-\left(\begin{array}{c}n+k-3\\ k-2\end{array}\right)
$$
if $k\ge 2$, whereas $m=n$ if $k=1$. 
We claim that there exists an eigenfunction 
\begin{equation}
\label{eq:eigenfunc}
\f(\sigma)=c_{1}H_{k,1}(\sigma)+\cdots+c_{m}H_{k,m}(\sigma)\quad(c_{1},\ldots,c_{m}\in\R)
\end{equation}
which vanishes on an equator of $\S^{n-1}$. More precisely, we can find $c_{1},\ldots,c_{m}\in\R$, not all zero, such that
\begin{equation}
\label{eq:algebraic-system}
c_{1}H_{k,1}(e_{j})+\cdots+c_{m}H_{k,m}(e_{j})=0\quad\forall j=1,\ldots,n-1
\end{equation}
where the $e_{j}$'s constitute the standard basis in $\R^{n}$. This is possible since $m>n-1$ as $n\ge 2$, via induction. Hence the algebraic system (\ref{eq:algebraic-system}) admits a nontrivial solution $c_{1},\ldots,c_{m}$. Thus the corresponding mapping $\f$ defined by (\ref{eq:eigenfunc}) restricted to 
$$
\S^{n-1}_{+}=\{\sigma\in\S^{n-1}~|~\sigma\cdot e_{n}>0\}
$$ 
is an eigenfunction in $H^{1}_{0}(\S^{n-1}_{+})$ relative to the eigenvalue $\lambda$. In fact $\f\ne 0$ since 
$$
2\|\f\|^{2}_{L^{2}(\S^{n-1}_{+})}=\|\f\|^{2}_{L^{2}(\S^{n-1})}=c_{1}^{2}+\cdots+c_{m}^{2}\ne 0.
$$
Therefore $\Lambda(\S^{n-1}_{+})$ contains $\Lambda(\S^{n-1})\setminus\{0\}$, and the proposition is completely proved.
\QED

\begin{Remark}
When $\alpha=0$ we obtain 
$$
\mu_{N}(\R^{n}_{+};0)=\left(\frac{n^{2}-4}{4}\right)^{2}.
$$
In particular $\mu_{N}(\R^{n}_{+};0)>0$ if and only if $n\ge 3$.
\end{Remark}

\subsection{Rellich inequality on cones in low dimension}
\label{SS:low}

For every $\theta\in(0,\pi)$ let $\Sigma_{\theta}\subset\S^{n-1}$ be a geodesic ball of radius $\theta$. Let us denote $\mathcal{C}_{\theta}=\mathcal{C}_{\Sigma_{\theta}}$. We investigate the dependence of the Rellich constant on $\theta$ for $\alpha=0$ in low dimensions. 

Firstly consider the dimension $n=2$, when $\gamma_{2,0}=-1$. From a direct computation of the spectrum of the Laplace-Beltrami operator on $\Sigma_{\theta}$ (see, e.g., \cite{CalMus}), it follows that 
$$
\mu_{N}(\mathcal{C}_{\theta};0)=\min\left\{\left(\frac{\pi^{2}}{4\theta^{2}}-1\right)^{2},\left(\frac{\pi^{2}}{\theta^{2}}-1\right)^{2}\right\}.
$$
Notice that $\theta^{*}:=\pi\sqrt{5/2}$ is a local maximum for the map $\theta\mapsto\mu_{N}(\mathcal{C}_{\theta};0)$. In addition, by (\ref{eq:monotone}),
$$
\mu_{D}(\R^{2}_{+};0)\ge \mu_{D}(\mathcal{C}_{\theta^{*}};0)>\mu_{N}(\mathcal{C}_{\theta^{*}};0)=\frac{9}{25}.
$$

When $n=3$ a similar phenomenon appears. In particular, there is exactly one value $\theta^{*}\in(\pi/2,\pi)$ such that
$-\gamma_{3,0}=3/4$ is the smallest eigenvalue of $-\Delta_{\sigma}$ on $\Sigma_{\theta^{*}}$. Thus $\mu_{N}(\mathcal{C}_{\theta^{*}};0)=0$ and $\mu_{N}(\mathcal{C}_{\theta};0)>0$ for $\theta<\theta^{*}$.

Finally we point out that in dimension $n=4$ it holds that $\mu_{N}(\mathcal{C}_{\Sigma};0)>0$ for all $\Sigma\ne\S^{3}$, since $\gamma_{4,0}=0$ and $\lambda_{\Sigma_{\theta}}>0$ for all $\theta\in(0,\pi)$.

\section{Inequalities with logarithmic weights}
\label{S:log}

In this Section we are concerned with inequalities, with sharp constants, involving the $L^{2}$ norm of $\Delta u$ with a weight $|x|^{\alpha}$, for mappings $u$ supported by cone-like domains. More precisely,
in this section we assume that
\begin{equation}
\label{eq:conelike}
\Omega=\mathcal{C}_{\Sigma}\cap\mathbb{B}^{n}=\{r\sigma~|~0<r<1,~\sigma\in\Sigma\}~,
\quad \textrm{or}\quad \Omega=\mathcal{C}_{\Sigma}\setminus\overline{\mathbb{B}^{n}}=\{r\sigma~|~r>1,~\sigma\in\Sigma\},
\end{equation}
where $\Sigma\subset\S^{n-1}$ is a domain of class $C^2$.
Notice that the case of the punctured ball
$\mathbb{B}^{n}\setminus\{0\}$ and the exterior domain 
$\R^{n}\setminus\overline{\mathbb{B}^{n}}$ are included by taking $\Sigma=\S^{n-1}$.

As in the previous Sections, we denote by $\lambda_{\Sigma}$ the first eigenvalue of the Laplace-Beltrami operator in $H^{1}_{0}(\Sigma)$ and we define $\gamma_{n,\alpha}$ and
$\overline\gamma_{n,\alpha}$ as in (\ref{eq:gammaNalpha}) and (\ref{eq:gammabar}),
respectively. We have the following result.

\begin{Theorem}
\label{T:log}
Let $\alpha\in\R$ and let $\Sigma$ be a domain of class $C^{2}$ in $\S^{n-1}$, with $n\ge 2$.
Let $\Omega$ be a domain as in (\ref{eq:conelike}).
Then:
\begin{itemize}
\item[(i)] \emph{[Navier case]}
For every $u\in C^{2}_{c}(\overline{\Omega}\setminus\{0\})$ it holds that
\begin{equation}
\label{eq:log-navier}
\int_{\Omega}
|x|^{\alpha}|\Delta u|^{2}
-\mu_{N}(\mathcal{C}_{\Sigma};\alpha)\int_{\Omega}|x|^{\alpha-4}|u|^2\ge\frac{\overline\gamma_{n,\alpha}+\lambda_{\Sigma}}{2}\int_{\Omega}|x|^{\alpha-4}|\log|x||^{-2}|u|^2.
\end{equation}
\item[(ii)] \emph{[Mixed Navier-Dirichlet case]}
For every $u\in C^{2}_{c}(\overline{\Omega}\setminus\{0\})$ with $\nabla u=0$ on $\S^{n-1}\cap\partial\Omega$, it holds that
\begin{equation}
\label{eq:log-navdir}
\begin{split}
\int_{\Omega}&
|x|^{\alpha}|\Delta u|^{2}-\mu_{N}(\mathcal{C}_{\Sigma};\alpha)\int_{\Omega}|x|^{\alpha-4}|u|^2\\
&\ge\frac{\overline\gamma_{n,\alpha}+\lambda_{\Sigma}}{2}\int_{\Omega}|x|^{\alpha-4}|\log|x||^{-2}|u|^2
+\frac{9}{16}\int_{\Omega}|x|^{\alpha-4}|\log|x||^{-4}|u|^2.
\end{split}
\end{equation}
\item[(iii)] \emph{[Dirichlet case]}
For every $u\in C^{2}_{c}(\Omega)$ it holds that
\begin{equation}
\label{eq:log-dirich}
\begin{split}
\int_{\Omega}&|x|^{\alpha}|\Delta u|^{2}
-\mu_{D}(\mathcal{C}_{\Sigma};\alpha)\int_{\Omega}|x|^{\alpha-4}|u|^2\\
&\ge\frac{\overline\gamma_{n,\alpha}+\lambda_{\Sigma}}{2}\int_{\Omega}|x|^{\alpha-4}|\log|x||^{-2}|u|^2
+\frac{9}{16}\int_{\Omega}|x|^{\alpha-4}|\log|x||^{-4}|u|^2.
\end{split}
\end{equation}
\end{itemize}
\end{Theorem}

\proof
We first handle the case $\Omega=\mathcal{C}_{\Sigma}\cap\mathbb{B}^{n}$. 
Fix $u\in C^{2}_{c}(\overline{\Omega}\setminus\{0\})$. Let $\mathcal{Z}_{\Sigma}^{+}=\R_+\times\Sigma$ and $w=Tu$, where $T$ is the Emden-Fowler transform defined in (\ref{eq:uw}). Arguing as in Section \ref{S:EF}, one can see that
$w\in H^{2}\cap H^{1}_{0}(\mathcal{Z}_{\Sigma}^{+}):=H^{2}(\mathcal{Z}_{\Sigma}^{+})\cap H^1_0(\mathcal{Z}_{\Sigma}^{+})$. In particular
\begin{equation}
\label{eq:LR1}
\int_{\Omega}|x|^{\alpha}|\Delta u|^{2}=\int_{\mathcal{Z}_{\Sigma}^{+}}|Lw|^{2}+\int_{\mathcal{Z}_{\Sigma}^{+}}|w_{ss}|^{2}+2\int_{\mathcal{Z}_{\Sigma}^{+}}|\nabla_{\!\sigma}w_{s}|^{2}+2\overline\gamma_{n,\alpha}\int_{\mathcal{Z}_{\Sigma}^{+}}|w_{s}|^{2}
\end{equation}
with $L$ defined as in (\ref{eq:Lw-cone}). Each term in the right hand side of (\ref{eq:LR1}) can be estimated according to the behaviour of $\nabla u$ on $\partial\Omega$. First of all, observe that for every $s>0$ the mapping $w(s,\cdot)$ belongs to $H^{2}\cap H^{1}_{0}(\Sigma)$. Hence we can apply Proposition \ref{P:auxiliary0} and Theorem \ref{T:Rellich-cone} to estimate
\begin{equation}
\label{eq:LR2}
\int_{\mathcal{Z}_{\Sigma}^{+}}|Lw|^{2}\ge{m}_{N}(\Sigma;\gamma_{n,\alpha})\int_{\mathcal{Z}_{\Sigma}^{+}}|w|^{2}=\mu_{N}(\mathcal{C}_{\Sigma};\alpha)\int_{\Omega}|x|^{\alpha-4}|u|^{2}.
\end{equation}
Similarly, if $u\in C^{2}_{c}(\Omega)$ then for a.e. $s>0$ the mapping $w(s,\cdot)$ belongs to $H^{2}_{0}(\Sigma)$ and in this case we obtain that
\begin{equation}
\label{eq:LR2zero}
\int_{\mathcal{Z}_{\Sigma}^{+}}|Lw|^{2}\ge{m}_{D}(\Sigma;\gamma_{n,\alpha})\int_{\mathcal{Z}_{\Sigma}^{+}}|w|^{2}=\mu_{D}(\mathcal{C}_{\Sigma};\alpha)\int_{\Omega}|x|^{\alpha-4}|u|^{2}.
\end{equation}
Notice also that for every $\sigma\in\Sigma$ the mapping $w(\cdot,\sigma)$ belongs to $H^{1}_{0}(\R_{+})$. Then we can exploit Hardy's inequality for functions of one variable
in order to estimate
\begin{equation}
\label{eq:LR3}
\int_{\mathcal{Z}_{\Sigma}^{+}}|w_{s}|^{2}\ge\frac{1}{4}\int_{\mathcal{Z}_{\Sigma}^{+}}s^{-2}|w|^{2}=\frac{1}{4}\int_{\Omega}|x|^{\alpha-4}|\log|x||^{-2}|u|^{2}.
\end{equation}
Moreover, since $w=0$ on $\R_{+}\times\partial\Sigma$, one has that $w_{s}(s,\cdot)\in H^{1}_{0}(\Sigma)$ for every $s>0$ and then
\begin{equation}
\label{eq:LR5}
\int_{\mathcal{Z}_{\Sigma}^{+}}|\nabla_{\!\sigma}w_{s}|^{2}\ge\lambda_{\Sigma}\int_{\mathcal{Z}_{\Sigma}^{+}}|w_{s}|^{2}.
\end{equation}
Then (\ref{eq:log-navier}) follows from (\ref{eq:LR2}), (\ref{eq:LR3}) and (\ref{eq:LR5}). If $\nabla u=0$ on $\S^{n-1}\cap\partial\Omega$ then for every $\sigma\in\Sigma$ the mapping $w(\cdot,\sigma)$ belongs to $H^{2}_{0}(\R_{+})$, that is, $w_{s}(\cdot,\sigma)\in H^{1}_{0}(\R_{+})$. Therefore, applying twice the Hardy inequality, we can estimate also
\begin{equation}
\label{eq:LR4}
\int_{\mathcal{Z}_{\Sigma}^{+}}|w_{ss}|^{2}\ge
\frac{9}{16}\int_{\mathcal{Z}_{\Sigma}^{+}}s^{-4}|w|^{2}=\frac{9}{16}\int_{\Omega}|x|^{\alpha-4}|\log|x||^{-4}|u|^{2}.
\end{equation}
Hence (\ref{eq:log-navdir}) follows from (\ref{eq:LR2}), (\ref{eq:LR3}), (\ref{eq:LR4}) and (\ref{eq:LR5}) whereas (\ref{eq:log-dirich}) follows from (\ref{eq:LR2zero}), (\ref{eq:LR3}), (\ref{eq:LR4}) and (\ref{eq:LR5}). 
Thus the theorem is proved when $\Omega$ is the intersection of a cone with the unit ball.

If $\Omega=\mathcal{C}_{\Sigma}\setminus\overline{\mathbb{B}^{n}}$ 
the proof is the same with the only difference that $\mathcal{Z}_{\Sigma}^{+}$ has to be changed into $\mathcal{Z}_{\Sigma}^{-}=\R_{-}\times\Sigma$. All the integrals in $w$ are invariant with respect to $s\mapsto -s$.
\QED

\begin{Remark}
Taking into account the role of the dilations in $\R^{n}$, Theorem \ref{T:log} can be suitably extended to any cone-like domain of the form 
$$\{r\sigma~|~0<r<R,~\sigma\in\Sigma\}
\quad\textrm{or}\quad
\{r\sigma~|~r>R,~\sigma\in\Sigma\}$$ with $R>0$ fixed. 
In fact the result $(iii)$ can be suitably extended to any bounded domain $\Omega$ with $0\in\partial\Omega$ or to any exterior domain $\Omega$ with $0\not\in\overline\Omega$, with no regularity assumption on $\partial\Omega$.  
\end{Remark}

In the next corollaries we point out the explicit constants in case $\alpha=0$, under Navier
and Dirichler boundary conditions. For the convenience
of the reader we distinguish the case $n=2$ from the higher dimensional one.

\begin{Corollary}
Let $\Omega=\mathbb{B}^{2}\setminus\{0\}$ or $\Omega=\R^{2}\setminus\overline{\mathbb{B}^{2}}$.
Then the following inequalities hold:
\begin{gather*}
\int_{\Omega}|\Delta u|^{2}\ge\frac{1}{2} \int_{\Omega}|x|^{-4}|\log|x||^{-2}|u|^2
\quad\textit{for any $u\in C^{2}_{c}(\overline\Omega\setminus\{0\})$,}\cr
\int_{\Omega}|\Delta u|^{2}\ge\frac{1}{2}\int_{\Omega}|x|^{-4}|\log|x||^{-2}|u|^2+
\frac{9}{16}\int_{\Omega}|x|^{-4}|\log|x||^{-4}|u|^2~\textit{for any $u\in C^{2}_{c}(\Omega)$.}
\end{gather*}
\end{Corollary}

\begin{Corollary}
Let $n\ge 3$. If $\Omega=\mathbb{B}^{n}\setminus\{0\}$ or $\Omega=\R^{n}\setminus\overline{\mathbb{B}^{n}}$ then the following inequalities hold:
$$
\int_{\Omega}|\Delta u|^{2}-\left(\frac{n(n-4)}{4}\right)^2\int_{\Omega}|x|^{-4}|u|^2\ge\frac{n^{2}-4n+8}{8}\int_{\Omega}|x|^{-4}|\log|x||^{-2}|u|^2
$$
for any $u\in C^{2}_{c}(\overline\Omega\setminus\{0\})$, and
\begin{eqnarray*}
\int_{\Omega}|\Delta u|^{2}&-&\left(\frac{n(n-4)}{4}\right)^2\int_{\Omega}|x|^{-4}|u|^2\\
&\ge&\frac{n^{2}-4n+8}{8}\int_{\Omega}|x|^{-4}|\log|x||^{-2}|u|^2+\frac{9}{16}\int_{\Omega}|x|^{-4}|\log|x||^{-4}|u|^2
\end{eqnarray*}
for any $u\in C^{2}_{c}(\Omega)$.
\end{Corollary}

In the next result we show that the constants appearing in the right hand side in (\ref{eq:log-navdir}) are sharp.

\begin{Theorem}
\label{T:log-sharp}
Let $\alpha$, $\Sigma$ and $\Omega$ as in Theorem \ref{T:log}, and 
Assume that
$$
\mathrm{dist}(-\gamma_{n,\alpha},\Lambda(\Sigma))=|\gamma_{n,\alpha}+\lambda_{\Sigma}|.
$$
If $A,B\in\R$ are such that 
\begin{equation}
\label{eq:Leray-Rellich-AB}
\begin{split}
\int_{\Omega}|x|^{\alpha}|\Delta u|^{2}-&\mu_{N}\mathcal{C}_{\Sigma};\alpha)\int_{\Omega}|x|^{\alpha-4}|u|^2\\
&\ge A\int_{\Omega}|x|^{\alpha-4}|\log|x||^{-2}|u|^2+B\int_{\Omega}|x|^{\alpha-4}|\log|x||^{-4}|u|^2
\end{split}
\end{equation}
for every $u\in C^{2}_{c}(\overline{\Omega}\setminus\{0\})$ with $\nabla u=0$ on $\S^{n-1}\cap\partial\Omega$, then 
$A\le(\overline\gamma_{n,\alpha}+\lambda_{\Sigma})/{2}$ and $B\le(3/4)^2$.
\end{Theorem}

\proof
Choose
$$
u(x)=|x|^{\frac{4-n-\alpha}{2}}v(-\log|x|)\varphi(\sigma)\quad(\sigma=x/|x|)
$$
where $\varphi$ is an eigenfunction corresponding to $\lambda_{\Sigma}$ and $v\in C^{2}_{c}(\R_{+})$. We compute
$$
\int_{\Omega}|x|^{\alpha}|\Delta u|^{2}=\int_{\Sigma}|\varphi|^{2}d\sigma\int_{0}^{\infty}\left(|v''|^{2}+2(\overline\gamma_{n,\alpha}+\lambda_{\Sigma})|v'|^{2}+(\gamma_{n,\alpha}+\lambda_{\Sigma})^{2}|v|^{2}\right)ds.
$$
Since 
$(\gamma_{n,\alpha}+\lambda_{\Sigma})^{2}=\mu_{N}(\mathcal{C}_{\Sigma};\alpha)$, 
(\ref{eq:Leray-Rellich-AB}) yields
\begin{equation}
\label{eq:LR-v}
\int_{0}^{\infty}|v''|^{2}ds+2(\overline\gamma_{n,\alpha}+\lambda_{\Sigma})\int_{0}^{\infty}|v'|^{2}ds\ge A\int_{0}^{\infty}s^{-2}|v|^2ds+B\int_{0}^{\infty}s^{-4}|w|^2ds
\end{equation}
for all $v\in C^{2}_{c}(\R_{+})$. As (\ref{eq:LR-v}) holds true for $w(s)=v(t s)$, for all $t>0$, we obtain
\begin{equation}
\label{eq:LR-vt}
t^{2}\int_{0}^{\infty}|v''|^{2}ds+2(\overline\gamma_{n,\alpha}+\lambda_{\Sigma})\int_{0}^{\infty}|v'|^{2}ds\ge A\int_{0}^{\infty}s^{-2}|v|^2ds+Bt^{2}\int_{0}^{\infty}s^{-4}|w|^2ds.
\end{equation}
Taking $t\to 0$ we get
$$
\frac{A}{2(\overline\gamma_{n,\alpha}+\lambda_{\Sigma})}\le
\inf_{\genfrac{}{}{0pt}{}{\scriptstyle{v\in C^{2}_{c}(\R_{+})}}{\scriptstyle{v\neq 0}}}
\frac{\displaystyle\int_{0}^{\infty}|v'|^{2}}{\displaystyle\int_{0}^{\infty}s^{-2}|v|^{2}}=\frac{1}{4}
$$
and thus we obtain the bound on $A$. Dividing (\ref{eq:LR-vt}) by $t^{2}$ and passing to the limit $t\to\infty$ we obtain
$$
B\le\inf_{\genfrac{}{}{0pt}{}{\scriptstyle{v\in C^{2}_{c}(\R_{+})}}{\scriptstyle{v\neq 0}}}\frac{\displaystyle\int_{0}^{\infty}|v''|^{2}}{\displaystyle\int_{0}^{\infty}s^{-4}|v|^{2}}=\frac{9}{16}.
$$
This completes the proof.\QED

Similar results on the optimality of the constants can be proved in the Navier and in the
Dirichlet case.

\label{References}

\noindent
\small
P. Caldiroli\\
Dipartimento di Matematica, Universit\`a di Torino\\
via Carlo Alberto, 10 -- 10123 Torino, Italy\\
\textit{E-mail address:} {paolo.caldiroli@unito.it}\\
\\
R. Musina\\
Dipartimento di Matematica ed Informatica, Universit\`a di Udine\\ via delle Scienze, 206 -- 33100 Udine, Italy\\
\textit{E-mail address:} {roberta.musina@uniud.it}


\begin{thebibliography}{XX}
\footnotesize


\bibitem{All}
{Allegretto, W.},
Nonoscillation theory of elliptic equations of order $2n$,
\textit{Pacific J. Math.} 
\textbf{64} (1976), 1--16.



\bibitem{Ben}
{Bennett, D.M.},
An extension of Rellich's inequality,
\textit{Proc. Amer. Math. Soc.} 
\textbf{106} (1989), 987--993.

\bibitem{BerFarFerGaz} 
{Berchio, E., Farina, A., Ferrero, A., Gazzola, F.}, 
{Existence and stability of entire solutions to a semilinear fourth order elliptic problem},
preprint (2011).

\bibitem{BerGauMaz}
{Berger, M., Gauduchon, P., Mazet, E.},
Le Spectre d'une Vari\'et\'e Riemannienne,
Lecture Notes in Mathematics \textbf{194}, Springer, 1971.

\bibitem{CKN}
{Caffarelli, L., Kohn, R., Nirenberg, L.},
{First Order Interpolation Inequalities with Weight},
\textit{Compositio Math.} \textbf{53} (1984), 259--275.

\bibitem{CalMus}
{Caldiroli, P., Musina, R.}, 
On a class of 2-dimensional singular elliptic problems, 
\textit{Proc. Roy. Soc. Edinburgh Sect. A} 
\textbf{131} (2001), 479--497.

\bibitem{CM2}
{Caldiroli, P., Musina, R.}, 
Caffarelli-Kohn-Nirenberg type inequalities for the weighted biharmonic operator and related semilinear problems. 
In progress.

\bibitem{DavHin}
{Davies, E.B., Hinz, A.M.},
Rellich inequalities in $L_{p}(\Omega)$,
\textit{Math. Z.} 
\textbf{227} (1998), 511--523. 

\bibitem{EvaLew}
{Evans, W.D., Lewis, R.T.},
On the Rellich inequality with magnetic potentials,
\textit{J. Math. Inequal.} 
\textbf{1} (2007), 473--490. 


\bibitem{GazGruMit} Gazzola, F., Grunau, H.-C., Mitidieri, E.,
Hardy inequalities with optimal constants and remainder terms,
\textit{Trans. Amer. Math. Soc.} 
\textbf{356} (2003), 2149--2168.

\bibitem{GazGruSwe} Gazzola, F., Grunau, H.-C., Sweers, G.,
Optimal Sobolev and Hardy--Rellich constants under Navier boundary conditions,
\textit{Ann. Mat. Pura Appl. (4)} 
\textbf{189} (2010), 475--486.

\bibitem{Mit} 
Mitidieri, E., 
A simple approach to Hardy's inequalities,
\textit{Math. Notes} 
\textbf{67} (2000), 479--486, 
translation from 
\textit{Mat. Zametki} 
\textbf{67} (2000), 563--572.
%

\bibitem{Rel54}
{Rellich, F.}, 
Halbbeschr\"ankte Differentialoperatoren h\"oherer Ordnung. 
In: J.C.H. Gerretsen, J. de Groot (Eds.): Proceedings of the International Congress of Mathematicians 1954, Volume III 
(pp. 243--250) Groningen: Noordhoff 1956.

\bibitem{Rel69}
{Rellich, F.},
{Perturbation theory of eigenvalue problems},
Gordon and Breach, New York, 1969.





\end{thebibliography}
\end{document}